\documentclass[12pt,twoside,bezier,reqno]{amsart}
\usepackage{epsfig}
\usepackage{amsmath,amsthm,amsfonts,amssymb,amscd,euscript}

\input{epsf}
\newcommand{\filebegin}{\begin{document}}
\newcommand{\fileend}{\end{document}}
\makeatother
\def\thefootnote{}
\newcommand{\lo}{\longrightarrow}
\newcommand{\NMM}{\hspace*{2mm}}

\renewcommand{\baselinestretch}{1.2}

\input{epsf}
\renewcommand{\baselinestretch}{1.1}
\def\n{\noindent}%

\numberwithin{equation}{section}
\def\mapdown#1{\Big\downarrow\rlap
{$\vcenter{\hbox{$\scriptstyle#1$}}$}}
\newtheorem{theorem}{Theorem}[section]
\newtheorem{lemma}[theorem]{Lemma}
\theoremstyle{definition}
\newtheorem{definition}[theorem]{Definition}
\newtheorem{example}[theorem]{Example}
\newtheorem{xca}[theorem]{Exercise}
\theoremstyle{remark}
\newtheorem{remark}[theorem]{Remark}
\numberwithin{equation}{section}
\theoremstyle{definition}
\newtheorem{thm}{Theorem}[section]
\newtheorem{cor}[thm]{Corollary}
\newtheorem{lem}[thm]{Lemma}
\newtheorem{prop}[thm]{Proposition}
\theoremstyle{definition}
\newtheorem{defn}[thm]{Definition}
\theoremstyle{remark}
\newtheorem{rem}[thm]{\bf Remark}
\theoremstyle{remark}
\numberwithin{equation}{section}
\newtheorem{exam}[thm]{\bf Example}
\newcommand{\BX}{\mathbf{B}(X)}
\newcommand{\A}{\mathcal{A}}
\newcommand{\PP}{\mathcal{P}}
\newcommand{\f}{\mathcal{F}}
\newcommand{\K}{\mathcal{K}}
\newcommand{\V}{\mathcal{V}}
\newcommand{\G}{\mathcal{G}}
\newcommand{\T}{\mathcal{T}}
\newcommand{\R}{\mathcal{R}}
\newcommand{\W}{\mathcal{W}}
\newcommand{\D}{\mathcal{D}}
\newcommand{\LL}{\mathcal{L}}
\newcommand{\Z}{\mathcal{Z}}
\newcommand{\X}{\mathcal{X}}
\newcommand{\N}{\mathbb{N}}
\newcommand{\h}{\mathcal{H}}
\newcommand{\C}{\mathcal{C}}
\newcommand{\B}{\mathcal{B}}
\newcommand{\db}{\overline{d}}
\newcommand{\du}{\underline{d}}
\newcommand{\I}{\mathcal{I}}
\newcommand{\J}{\mathcal{J}}
\newcommand{\ub}{\overline{U}}
\newcommand{\fs}{\mathcal{P}_{f}}
\newcommand{\p}{\widetilde{p}}
\newcommand{\q}{\widetilde{q}}
\newcommand{\rr}{\widetilde{r}}
\newcommand{\LM}{Lmc(S)}
\newcommand{\LU}{\mathcal{LUC}}

\begin{document}


\vspace*{2cm}
\begin{center}
{\bf\large Central set Theorem near zero}
 \\[0.5cm]
{E. Bayatmanesh\\
M. Akbari Tootkaboni\\[2mm]
Department of Mathematics, Faculty of Basic Science,\\
 Shahed University,\\
  Tehran-Iran.\\[2mm]
{\tt E-mail: bayatmanesh.e$@$gmail.com}\\[2mm]
{\tt E-mail: akbari$@$shahed.ac.ir}}\\[2mm]
{A. Bagheri Sales\\[2mm]
 Department of Mathematics, Faculty of science,\\
Qom University,\\
Qom-Iran.\\[2mm]
{\tt E-mail: alireza{\underline{ }}bagheri{\underline{ }}salec@yahoo.com}}\\[2mm]
\end{center}
\vspace*{0.5cm}
\begin{quotation}
\noindent
{\footnotesize
{\sc Abstract.}
In this paper, we introduce notions of
$J$-set near zero and $C$-set near zero for a dense subsemigroup of $((0,+\infty),+)$
 and obtain some results for them.
Also we derive the Central Sets Theorem near
zero.}
\end{quotation}

\n \textbf{2000 Mathematics subject classification: } Primary: 54D80, 22A15, Secondary: 22A20.\\

\markboth
{E. Bayatmanesh, M. A. Tootkaboni}
 {Central set Theorem near zero}



\section{Introduction}

Let $(S,+)$ be a discrete semigroup. The collection of all ultrafilters on $S$ is called the Stone-$\check{C}$ech compactification
of $S$ and denoted by $\beta S$. For $A\subseteq S$, define $\overline{A}=\{p\in\beta S:A\in p\}$, then $\{\overline{A}:A\subseteq S\}$ is
a basis for the open sets( also for the closed sets) of $\beta S$. There is a unique extension of the operation to $\beta S$ making $(\beta S,+)$ a right
topological semigroup ( i.e. for each $p\in \beta S$, the right translation $\rho_p$ is continuous where $\rho_p(q)=q+p$) and also
for each $x\in S$, the left translation $\lambda_x$ is continuous where $\lambda_x(q)=x+q$. The principal
ultrafilters being identified with the points of $S$ and $S$ is a dense subset of $\beta S$. Given $p,q\in\beta S$ and $A\subseteq S$, we have $A\in p+q$ if and only if
$\{x\in S:-x+A\in q\}\in p$, where $-x+A=\{y\in S:x+y\in A\}$.

A nonempty subset $L$ of a semigroup $(S,+)$ is called a left ideal of $S$ if
$S + L \subseteq L$, a right ideal if $L + S \subseteq L$, and a two sided ideal (or simply an
ideal ) if it is both a left and right ideal. A minimal left ideal is a left ideal
that does not contain any proper left ideal. Similarly, we can define minimal
right ideal and smallest ideal.

Any compact Hausdorff right topological semigroup $(S,+)$ has a smallest
two sided ideal, denoted by $K(S)$,which is  the union of all minimal left ideals.  $K(S)$ is the
union of all minimal right ideals, as well. Given a minimal left ideal $L$ and a minimal right ideal $R$, $L\cap R$ is a group,
and in particular contains an idempotent. An idempotent in $K(S)$ is called
a minimal idempotent. An idempotent is minimal if and only if it is a member of the smallest ideal. For more details see \cite{hinbook}.

 For $A\subseteq S$, and $p\in \beta S$, we define
$A^{*}(p)=\{s\in A: -s+A\in p\}$.
\begin{lemma}\label{lemma11}
Let $(S,+)$ be a semigroup, let $p+p=p\in \beta S$, and let $A\in p$. For each $s\in A^{*}(p)$, $-s+A^{*}(p)\in p$.
\end{lemma}
\begin{proof}
\cite {hinbook}, Lemma 4.14.
\end{proof}
Now we review the definition of partition regularity and a theorem that connects it with ultrafilters.
   \begin{defn}
  Let $\R$ be a nonempty set of subsets of $S$.  $\R$ is partition regular if and only if
   whenever $\f$ is a finite set of $\mathcal{P}(S)$ and $\bigcup\f\in\R$, there exist $A\in\f$ and
    $B\in\R$ such that $B\subseteq A$.
     \end{defn}

\begin{thm}\label{theorem0}
Let $\R\subseteq \mathcal{P}(S)$  be a nonempty set
and assume $\emptyset\notin\R$.
 Let $$\R^\uparrow=\{B\in \mathcal{P}(S):A\subseteq B\mbox{ for some }A\in\R\}.$$
 Then (a), (b) and (c) are equivalence.

 (a) $\R$ is partition regular.

 (b) Whenever $\A\subseteq \mathcal{P}(S)$ has the property that every finite nonempty  subfamily of
 $\A$ has an intersection which is in $\R^\uparrow$, there is
  $ \mathcal{U} \in \beta S_d$ such that $\A\subseteq \mathcal{U} \subseteq \R^\uparrow$.

 (c) Whenever $A\in\R$, there is  $ \mathcal{U}\in\beta S_d$ such that $A\in \mathcal{U} \subseteq \R^\uparrow$.
 \end{thm}
{\flushleft\bf Proof.} \cite[Theorem 3.11]{hinbook}.$\hfill\square$

\begin{defn}
Let $(S,.)$ be a  discrete semigroup, let $A\subseteq S$. Then $A$ is a central set if and only if there exists an idempotent $p$
in the smallest  ideal of $\beta S$ with $A\in p$.
\end{defn}
In this paper, the collection of all nonempty subsets of $S$ is denoted by $P_f(S)$ and
$ P(S) $ is the set of all subsets of $S$. Central subsets of a discrete semigroup $S$
have very strong combinatorial properties which are  consequences of the central set theorem.
There is an elementary description of central sets which was showed in \cite{hinbook}.
\begin{thm}[Central Set Theorem \cite {Di-hindman}]
Let $(S,+)$ be a commutative  subsemigroup  and let $\T={}^{\Bbb N}S$, the set of sequences in $S$. Let $A$ be a central subset of $S$.
There exist functions $\alpha:P_f(\T)\rightarrow S$ and $H:P_f(\T)\rightarrow P_f(\N)$ such that

(1) if $F,G\in P_f(\T)$ and $F\subsetneq G$, then $maxH(F)<min H(G)$, and

(2) whenever $m\in \mathbb{N}$, $G_1,\cdots, G_m\in P_f(\T)$, $G_1\subsetneq G_2\subsetneq\cdots \subsetneq G_m$, and
for each $i\in \{1,2,\cdots , m\}$, $f_i\in G_i$, one has
\[
\sum_{i=1}^m\big(\alpha(G_i)+\sum_{t\in H(G_i)}f_i(t)\big)\in A.
\]
\end{thm}
We define a set to be  a $C-$set if and only if it satisfies the conclusion of  Central Set Theorem.
\begin{defn}\label{defn1}
Let $S$ be a commutative subsemigroup  and let $A\subseteq S$ and let $\T=^\N S$. The set $A$ is a $C-$set  if and only if
there exist functions $\alpha:P_f(\T)\rightarrow S$ and $H:P_f(\T)\rightarrow P_f(\N)$ such that

(1) if $F,G\in P_f(\T)$ and $F\subsetneq G$, then $maxH(F)<min H(G)$, and

(2) whenever $m\in \mathbb{N}$, $G_1,\cdots, G_m\in P_f(\T)$, $G_1\subsetneq G_2\subsetneq\cdots \subsetneq G_m$, and
for each $i\in \{1,2,\cdots , m\}$, $f_i\in G_i$, one has
\[
\sum_{i=1}^m\big(\alpha(G_i)+\sum_{t\in H(G_i)}f_i(t)\big)\in A.
\]
\end{defn}
The central sets are important, and algebraically   are easy to work with.
However, from a combinatorial viewpoint, $C-$sets are the objects that matter. In $(\N,+)$ and in
many other semigroups, they are the objects that contain solutions to partition regular systems of
homogeneous equations as well as the other myriads of properties that are consequences of  Central
Sets Theorem, see e.g. \cite{hin2}.  In \cite{Di-hindman} authors obtained a simple characterization of $C-$sets in an arbitrary discrete semigroup.

\begin{defn}
Let $(S,+)$ be a commutative subsemigroup and let $\T=^\N S$. Let $A\subseteq S$  is a
$J-$set  if and only if whenever $F\in P_f(\T)$, there exist $a\in S$ and $H\in P_f(\N)$
such that for each $f\in F$, $a+\sum_{t\in H}f(t)\in A$.
\end{defn}
If $S$ is noncommutative, then the definition of $J-$sets is somewhat more complicated, but still much simpler than the statement of the
noncommutative Central Sets Theorem. One of the reasons that $J-$sets are of
interest is that, if $S$ is a discrete commutative semigroup,
then every subset of $S$ with positive upper density is a $J-$set, see \cite{hin2}.

In \cite{Hin-Lead}, N. Hindman and I. Leader investigated  some concepts near zero. The set $0^+$ of all non principal ultrafilters on $S=((0,\infty),+)$ that is convergent to $0$  is a semigroup under the restriction of the usual $'+'$ on $ \beta S$, the Stone-$\check{C}$Cech compactification of the discrete semigroup $S=((0,\infty),+)$. Further, in \cite{Di-hindman2} and \cite{hind-de}, De, Hindman and Strauss worked out
 some more and much more elaborated applications of central sets. One such property is that in any finite partition of a central set at least one
 cell of the partition is central set \cite {Furesten}. In \cite{Di-hindman2} the authors used the algebraic structure
 of $0^+(S)$ in their investigation of image partition regularity near $0$ of finite and infinite matrices.

In \cite{D-paul} the algebraic structure of $0^+(\mathbb{R})$  is used to investigate image partition regularity  of matrices with real entries form $\mathbb{R}$.
Central sets near zero were introduced by N. Hindman and I. Learder  in \cite{Hin-Lead} as central sets, central sets near zero  enjoy rich combinatorial structure ,too.
Central sets are ideal objects for Ramsey theoretic. The central sets theorem was first introduced by H. Furstenberg ( see \cite {Furesten}) for
the semigroup $\N$ and considering sequence in $\mathbb{Z}$. The most general version of Central Sets Theorem is available in \cite{Di-hindman}.

 We have been considering semigroups which are dense in $((0,\infty),+)$ with natural topology.
When passing to the Stone-$\check{C}$ech compactification of such a semigroup S, we deal with
$S_d$, which is the set $S$ with the discrete topology.
\begin{defn}
Let $S$ be a dense subset of $((0,\infty),+)$. Then
\[
0^+(S)=\{p\in\beta S_d:(\forall\epsilon>0)\,\,(0,\epsilon)\cap S\in p\}.
\]
\end{defn}
By Lemma 2.5 in \cite{Hin-Lead}, $0^+(S)$ is a compact right topological
subsemigroup of $(\beta S_d,+)$, and $0^+(S)\cap K(\beta S_d)=\emptyset$. Since $0^+(S)$
is compact right topological semigroup, so $0+(S)$ contains minimal idempotents.

In section 2 we difine notions of thick near zero, syndetic near zero and  piecewise syndetic near zero and  
we obtain a main theorem about piecewise syndetic near zero.

In section 3,  we introduce Central Set Theorem near zero,  also 
we define  $C$-set near zero, $J$-set near zero and  derive some results for dense subsemigroups of $((0,\infty),+)$.
We will be considering semigroups which are dense  with respect to the usual topology on $((0,\infty),+)$.

\section{Additive properties near zero}
 In this section, we define  the notions of  thick near zero, syndetic near zero and piecewise
 syndetic near zero,  which   in turn    extend  the notions  of  thick,  syndetic
 and piecewise  syndetic. Recall that $A \subseteq S$   is thick if and only
  if for every $F\in\fs(S)$ there exists $x \in S$ such that $F+x \subseteq A$,    $A$
  is syndetic if and only if there exists $H\in\fs(S)$ such that $S =\bigcup_{t\in H} -t+A$ and $A $
  is piecewise syndetic if and only if there exists $H \in\fs(S)$ such that $\bigcup_{t\in H} -t+A$ is thick.
\vspace{.1cm}

 \begin{defn}
Let $S$ be a dense subsemigroup of $((0,\infty),+)$, and let $A \subseteq S$.

$(a)$  $A$ is thick near zero if and only if
$(\exists \epsilon>0) (\forall F\in (P_f(0,\epsilon)\cap S)) (\forall \delta>0) (\exists y\in (0,\delta)\cap S) (F+y\subseteq A)$.

  $(b)$ $A$ is syndetic near zero if and only if for any  $\epsilon>0$ there exist  $F\in\PP_f((0,\epsilon)\cap S)$ and
  $\delta>0$ such that $(0,\delta)\cap S \subseteq\bigcup_{t\in F} -t+A$.

   $(c)$ $A$ is piecewise syndetic near zero if and only if for all $\delta>0$ there exists $F\in\PP_f((0,\delta)\cap S)$
   such that $\bigcup_{t\in F}-t+A$ be thick near zero.
   
   $(d)$  $A$ is a central set near zero if and only if there exists an idempotent $p$
in the smallest  ideal of $0^+(S)$ with $A\in p$.

\end{defn}
\begin{rem}
It is obvious that $A$ is thick near zero if and only if for some $\delta>0$
\[
\{-t+A:t\in (0,\delta)\cap S\}
\]
has the finite intersection property in $(0,\delta)\cap S$, i.e. for each $t_1,\cdots ,t_n\in (0,\delta)\cap S$
there exists $x\in (0,\delta)\cap S$ such that $x\in \bigcap_{i=1}^n-t_i+A$.

\end{rem}
In this paper, the minimal ideal in $0^+(S)$ is denoted by $K$.
\begin{thm}
Let $ S $ be a dense subsemigroup of $((0,+\infty),+) $ and let $ p \in 0^+(S) $. The following statements are equivalent.

(a) $ p \in K $.

 (b) For all $ A \in p $, $ \left\lbrace x \in S : -x+A \in p\right\rbrace  $ is syndetic near $ 0 $.

 (c) For all $ r \in 0^+(S) $, $ p \in 0^+(S) +r+p $.
\end{thm}
\begin{proof}
See \cite{Hin-Lead} or see Theorem 3.4 in \cite{Ak-va}.
\end{proof}

\begin{thm}\label{theorem01}
Let $A\subseteq S$. Then $K\cap \overline{A}\neq\emptyset$ if and only if $A$ is piecewise syndetic near $0$.
\end{thm}
\begin{proof}
Necessity. Pick $p\in K\cap \overline{A}$ and let $B=\{x\in S:-x+A\in p\}$. By Theorem 2.3, $B$ is syndetic
near $\mu$. So for every $\varepsilon>0$ there exist $F\in\PP_f((0,\varepsilon)\cap S)$ and
$\delta>0$ such that $(0,\delta)\cap S\subseteq \bigcup_{t\in F}-t+B$. So for each $x\in (0,\delta)\cap S$,
there exists $t\in F$ such that $x\in -t+B$, and so $-(x+t)+A\in p$. Thus $-x+(\bigcup_{t\in F}-t+A)\in p$
and since $(0,\delta)\cap S\in p$ for each $\delta>0$, hence $\{-x+(\bigcup_{t\in F}-t+A):x\in (0,\delta)\cap S\}$ has the finite intersection property
in $(0,\delta)\cap S$. So by Remark 2.2, $A$ is piecewise syndetic near zero.

Sufficiency. Let $A$ be piecewise syndetic near zero. So for each $n\in \mathbb{N}$ there exists $F_n\in\PP_f((0,\frac{1}{n})\cap S)$ such that
for some $\epsilon_n>0$ and for each $G_n\in\PP_f((0,\epsilon_n)\cap S)$, for every $\delta_n$ there exists $y_n\in \PP_f((0,\delta_n)\cap S)$ such that $G_n+y_n\subseteq \bigcup_{t\in F_n}(-t+A)$. Pick $\delta_n=min({\dfrac{1}{n},\epsilon_n})$
For $G\in\PP_f(S)$  and $\mu>0$, let
\[
C(G,\mu)=\{x\in (0,\mu):\forall n\in \mathbb{N}, G\cap (0,\delta_n)+x\subseteq\bigcup_{t\in F_n}(-t+A)\}.
\]
It is obvious that each $C(G,\mu)\neq\emptyset$. Also,
\[
C(G_1\cup G_2,min\{\mu_1,\mu_2\})\subseteq C(G_1,\mu_1)\cap C(G_2,\mu_2),
\]
for each $G_1,G_2\in\PP_f(S)$ and for each $\mu_1,\mu>0$. Therefore
\[
\{C(G,\mu):G\in \PP_f(S)\mbox{ and }\mu>0\}
\]
 has the finite intersection property, so pick
$p\in \beta S_d $ with
\[
\{C(G,\mu):G\in \PP_f(S)\mbox{ and }\mu>0\}\subseteq p.
\]
Since $C(G,\mu)\subseteq (0,\mu)\cap S$, so $p\in 0^+(S)$.

Now we claim that for each $n\in \mathbb{N}$, $0^+(S)+p\subseteq cl_{\beta S_d}\bigcup_{t\in F_n}(-t+A)$, so let
$n\in \mathbb{N}$ and let $q\in 0^+(S)$. To show that $\bigcup_{t\in F_n}(-t+A)\in q+p$, we show that
\[
(0,\delta_n)\cap S\subseteq\{y>0:-y+\bigcup_{t\in F_n}(-t+A)\in p\}.
\]
So let $y\in (0,\delta_n)\cap S$. Then $C(\{y\},\delta_n)\in p$ and $C(\{y\},\delta_n)\subseteq -y+\bigcup_{t\in F_n}(-t+A)$.
Now pick $r\in (0^+(S)+p)\cap K$( since $0^+(S)+p$ is a left ideal of $0^+(S)$). Given $n\in \mathbb{N}$, $\bigcup_{t\in F_n}(-t+A)\in r$
so pick $t_n\in F_n$ such that $-t_n+A\in r$. Since for each $n\in\mathbb{N}$, $t_n\in F_n\subseteq (0,\frac{1}{n})\cap S$ so $lim_{n\rightarrow\infty}t_n=0$.
 Now pick $q\in 0^+(S)\cap cl_{\beta S_d}\{t_n:n\in \mathbb{N}\}$. Then $q+r\in K$ and
$\{t_n:n\in \mathbb{N}\}\subseteq \{t\in S: -t+A\in r\}$ so $A\in q+r$.
\end{proof}

\begin{cor}
Let $(S,+)$ be a dense subsemigroup of $(0,\infty)$,  piecewise syndeticity  near zero is partition regular.
\end{cor}
\begin{proof}
It is obvuous.
\end{proof}

\section{Central Sets}
Central subsets of a discrete semigroups have very strong combinatorial properties which are a consequence of the Central Sets Theorem.
\begin{defn}
Let $S$ be a dense subsemigroup of $\big((0,\infty),+\big)$. We say that
$f:\mathbb{N}\rightarrow S$ is near zero if $inf f(\N)=0$. The collection of all
functions that is near zero is denoted by $\T_0$.
\end{defn}
\begin{defn}
Let $S$ be a dense subsemigoup of $((0,\infty),+)$ and let $A\subseteq S$. Then $A$ is a
$J-$set near zero if and only if whenever $F\in P_f(\T_0)$ and $\delta>0$, there exist $a\in S\cap (0,\delta)$ and $H\in P_f(\N)$
such that for each $f\in F$, $a+\sum_{t\in H}f(t)\in A$.
\end{defn}
Of course, we can say that $A\subseteq S$ is a $J-$set near zero if and only if
for each $F\in P_f(\T_0)$ and for each $\delta>0$, there exist $a\in S$ and $H\in P_f(\N)$ such that
$a+\sum_{t\in H}f(t)\in A\cap (0,\delta)$ for each $f\in F$, i.e. for each $\delta>0$, $A\cap (0,\delta)$ is a $J-$set. It is obvious that
every $J-$set near zero respect to this definition is a $J-$set near zero by Definition 5.2. So we focuse on Definition 5.2.
\begin{lem}\label{lemma1}
Let $S$ be a dense subsemigroup of $\big((0,\infty),+\big)$ and let $A\subseteq S$ be a $J$-set near zero. Whenever $m\in\N$ and
$F\in P_f(\T_0)$ and $\delta>0$, there exist $a\in S\cap (0,\delta)$ and $H\in P_f(\N)$ such that $minH>m$ and for each $f\in F$, $a+\sum_{t\in H}f(t)\in A$.
\end{lem}
\begin{proof}
See Lemma 14.8.2 in \cite {hinbook}.
\end{proof}
\begin{thm}\label{theorem1}
Let $S$ be a dense subsemigroup of $\big((0,\infty),+\big)$ and let $A$ be a subset of $S$.
If $A$ is a piecewise syndetic near zero, then $A$ is a $J$-set near zero.
\end{thm}
\begin{proof}
Let $F\in P_f(\T_0)$, let $l=|F|$, and enumerate $F$ as $\{f_1,\cdots, f_l\}$. Let $Y=\prod_{t=1}^l0^+$.
Then by Theorem 2.22 in \cite{hinbook}, $Y$ is a compact right topological semigroup and if $s\in\prod_{t=1}^lS$,
then $\lambda_s$ is continuous. For $i\in\N$ and $\delta>0$, let

\begin{multline*}
  I_{i,\delta}=\big\{\big(a+\sum_{t\in H}f_1(t),\cdots,a+\sum_{t\in H}f_l(t)\big):a\in S\cap(0,\delta),\\
   H\in P_f(\N),\mbox{ and }minH>i\big\}\\
\end{multline*}
and let $E_{i,\delta}=I_{i,\delta}\cup\{(a,\cdots,a):a\in S\cap (0,\delta)\}$.

Let $E=\bigcap_{i\in\N,\delta>0}\overline{E_{i,\delta}}$ and let $I=\bigcap_{i\in\N,\delta>0}\overline{I_{i,\delta}}$. It is obvious that
$E\subseteq 0^+$ and $I\subseteq 0^+$. We claim
that $E$ is a subsemigroup of $Y$ and $I$ is an ideal of $E$. To this end, let $p,q\in E$. We show that
$p+q\in E$ and if either $p\in I$ or $q\in I$, then $p+q\in I$. Pick $\delta>0$, then $U=cl_{\beta S_d}\big((0,\delta)\cap S\big)$ is an open neighborhood
of $p+q$ and let $i\in\N$. Since $\rho_q$ is continuous, pick a neighborhood $V$ of $p$ such that
$V+q\subseteq U$. Pick $x\in E_{i,\frac{\delta}{3}}\cap V$ with $x\in I_{i,\frac{\delta}{3}}$ if $p\in I$. If $x\in I_{i,\frac{\delta}{3}}$ so
that $x=(a+\sum_{t\in H}f_1(t),\cdots,a+\sum_{t\in H}f_l(t))$ for some $a\in S\cap (0,\frac{\delta}{3})$ and some $H\in P_f(\N)$
with $minH>i$, let $j=maxH$. Otherwise, let $j=i$. Since $\lambda_x$ is continuous, pick a neighborhood
$W$ of $q$ such that $x+W\subseteq U$. Pick $y\in E_{j,\frac{\delta}{3}}\cap W$ with $y\in I_{j,\frac{\delta}{3}}$ if $q\in I$.
Then $x+y\in E_{i,\delta}\cap U$ and if either $p\in I$ or $q\in I$, then $x+y\in I_{i,\delta}\cap U$.

By Theorem 2.23 in \cite {hinbook}, $K(Y)=\prod_{t=1}^lK(0^+)$. Pick by Theorem 2.4 some $p\in K(0^+)\cap \overline{A}$.
Then $\overline{p}=(p,\cdots,p)\in K(Y)$. We claim that $\overline{p}\in E$. To see this, let $U$ be a neighborhood of
$\overline{p}$, let $i\in\N$, and pick $C_1,\cdots,C_l\in p$ such that $\prod_{t=1}^l\overline{C_t}\subseteq U$. Pick
$a\in\bigcap_{t=1}^lC_t$. Then $\overline{a}=(a,\cdots,a)\in U\cap E_{i,\delta}$. Thus
$\overline{p}\in K(Y)\cap E$ and consequently $K(Y)\cap E\neq\emptyset$. Then by Theorem 1.65 in \cite{hinbook}, we have that
$K(E)=K(Y)\cap E$ and so $\overline{p}\in K(E)\subseteq I$. Then $I_{1,\delta}\cap \prod_{t=1}^l\overline{A}\neq\emptyset$ for each
$\delta>0$, so pick $z\in I_{1,\delta}\cap\prod_{t=1}^l \overline{A}$ and pick $a\in S\cap (0,\delta)$ and $H\in P_f(\N)$ such that
\[
z=\big(a+\sum_{t\in H}f_1(t),\cdots,a+\sum_{t\in H}f_l(t)\big).
\]
\end{proof}
\begin{thm}[Central set Theorem near zero]
Let $S$ be a dense subsemigroup of $\big((0,\infty),+\big)$. Let $A$ be a central subset of $S$ near zero. Then for each $\delta\in (0,1)$,
there exist functions $\alpha_\delta:P_f(\T_0)\rightarrow S$ and $H_\delta:P_f(\T_0)\rightarrow P_f(\N)$ such that \\
(1) $\alpha_\delta(F)<\delta$ for each $F\in P_f(\T_0)$,\\
(2) if $F,G\in P_f(\T_0)$ and $F\subset G$, then $maxH_\delta(F)<min H_\delta(G)$ and \\
(3) whenever $m\in \mathbb{N}$, $G_1,\cdots, G_m\in P_f(\T_0)$, $G_1\subset G_2\subset\cdots \subset G_m$, and
for each $i\in \{1,2,\cdots , m\}$, $f_i\in G_i$, one has
\[
\sum_{i=1}^m\big(\alpha_\delta(G_i)+\sum_{t\in H_\delta(G_i)}f_i(t)\big)\in A.
\]
\end{thm}
\begin{proof}
Pick a minimal idempotent $p$ of $0^+$ such that $A\in p$. Let $A^*=\{x\in A:-x+A\in p\}$, so $A^*\in p$.
Also by Lemma 4.14 in \cite {hinbook}, if $x\in A^*$, then $-x+A^*\in p$.

We define $\alpha_\delta(F)\in S$ and $H_\delta(F)\in P_f(\N)$ for $F\in P_f(\T_0)$
by induction on $|F|$ satisfying the following inductive
hypotheses:\\
(1) $\alpha_\delta(G)<\delta$ for each $G\in\T_0$,\\
(2) if $F,G\in P_f(\T_0)$ and $F\subset G$, then $maxH_\delta(F)<min H_\delta(G)$ and \\
(3) whenever $m\in \mathbb{N}$, $G_1,\cdots, G_m\in P_f(\T_0)$, $G_1\subset G_2\subset\cdots \subset G_m$, and
for each $i\in \{1,2,\cdots , m\}$, $f_i\in G_i$, one has
\[
\sum_{i=1}^m\big(\alpha_\delta(G_i)+\sum_{t\in H_\delta(G_i)}f_i(t)\big)\in A^*.
\]

Assume that $F=\{f\}$. Since $A^*$ is piecewise syndetic near zero, pick by
Theorem \ref{theorem1}, for $\delta>0$, $a\in S\cap (0,\delta)$ and $L\in P_f(\N)$ such that
$a+\sum_{t\in L}f(t)\in A^*$. Let $\alpha_\delta(\{f\})=a$ and
$H_\delta(\{f\})=L$.

Let $|F|>1$, $\alpha_\delta(G)$ and $H_\delta(G)$ have been defined for all proper subsets $G$ of $F$ and for each $\delta>0$.
Pick $\delta>0$, and let
\[
K_\delta=\bigcup\{H_\delta(G):G\mbox{ is a non-empty proper subset of }F\}
\]
 and let $m=max K_\delta$. Let
\begin{eqnarray*}
M_\delta&=&\big\{\sum_{i=1}^n\big(\alpha_\delta(G_i)+\sum_{t\in H_\delta(G_i)}f_i(t)\big):n\in \mathbb{N},\emptyset\neq G_1\subset\cdots \subset G_n\subset F,\\
&and& \{f_i\}_{i=1}^n\in\prod_{i=1}^n G_i\big\}.
\end{eqnarray*}
Then $M_\delta$ is finite and by hypothesis (3), $M_\delta\subseteq A^*$. Let
$B=A^*\cap\bigcap_{x\in M_\delta}(-x+A^*)$. Then $B\in p$ so pick by
Theorem \ref{theorem1} and Lemma \ref{lemma1}, $a\in S\cap (0,\delta)$ and $L\in
P_f(\N)$ such that $a+\sum_{t\in L}f(t)\in B$ for each $f\in F$.
Let $\alpha_\delta(F)=a$ and $H_\delta(F)=L$.

The hypothesis (1) is obvious. Since $minL>m$, we have the hypothesis (2) is satisfied. To verify hypothesis (3), pick $\delta>0$ and
$n\in \mathbb{N}$, let $\emptyset\subset G_1\subset\cdots \subset G_n=F$, and let $\{f_i\}_{i=1}^n\in\prod_{i=1}^n G_i$.
If $n=1$, then $\alpha(G_1\prod\delta)+\sum_{t\in H(G_1)}f_1(t)=a+\sum_{t\in L}f_1(t)\in B\subseteq A^*$. So assume that $n>1$ and
let $y=\sum_{i=1}^{n-1}\big(\alpha_\delta(G_i)+\sum_{t\in H_\delta(G_i)}f_i(t)\big)$.
Then $y\in M_\delta$ so $a+\sum_{t\in L}f_1(t)\in B\subseteq(-y+A^*)$ and thus
$\sum_{i=1}^n\big(\alpha_\delta(G_i)+\sum_{t\in H_\delta(G_i)}f_i(t)\big)=y+a+\sum_{t\in L}f_1(t)\in A^*$ as required.
\end{proof}
\begin{defn}\label{defn1}
Let $S$ be a dense subsemigroup of $\big((0,\infty),+\big)$ and let $A\subseteq S$. We say $A$ is a $C-$set near zero if and only if
for each $\delta\in (0,1)$,
there exist functions $\alpha_\delta:P_f(\T_0)\rightarrow S$ and $H_\delta:P_f(\T_0)\rightarrow P_f(\N)$ such that \\
(1) $\alpha_\delta(F)<\delta$ for each $F\in P_f(\T_0)$,\\
(2) if $F,G\in P_f(\T_0)$ and $F\subset G$, then $maxH_\delta(F)<min H_\delta(G)$ and \\
(3) whenever $m\in \mathbb{N}$, $G_1,\cdots, G_m\in P_f(\T_0)$, $G_1\subset G_2\subset\cdots \subset G_m$, and
for each $i\in \{1,2,\cdots , m\}$, $f_i\in G_i$, one has
\[
\sum_{i=1}^m\big(\alpha_\delta(G_i)+\sum_{t\in H_\delta(G_i)}f_i(t)\big)\in A.
\]
\end{defn}
Let $\Phi$ is the set of all functions $f:\N\rightarrow\N$ for which $f(n)\leq n$ for each $n\in\N$.
\begin{thm}\label{theorem5.7}
Let $S$ be a dense subsemigroup of $\big((0,\infty),+\big)$ and let $A$ be a $C$-set near zero in $S$, and for each $l\in\N$,
let $\{y_{l,n}\}_{n\in\N}\in\T_0$. There exist a sequence $\{a_n\}_{n\in\N}$ in $S$
such that $a_n\rightarrow 0$ and a sequence $\{H_n\}_{n\in\N}$ in $P_f(\N)$ such that $maxH_n<minH_{n+1}$ for each $n\in\N$ and such that for each $f\in\Phi$
\[
FS\big(\{a_n+\sum_{t\in H_n}y_{f(n),t}\}_{n\in\N}\big)\subseteq A.
\]
In particular, the above conclusion applies if $A$ is a central set near zero in $S$.
\end{thm}
\begin{proof}
Pick $\alpha$ and $H$ as guaranteed by Definition \ref{defn1}. We may assume that the sequences $\{y_{l,n}\}_{n\in\N}$ are distinct.

For $n\in\N$, let $F_n=\{\{y_{l,t}\}_{t\in\N},\{y_{2,t}\}_{t\in\N},\cdots,\{y_{n,t}\}_{t\in\N}\}$
and let $a_n=\alpha_\frac{1}{n}(F_n)<\frac{1}{n}$ and $H_n=H_\frac{1}{n}(F_n)$. Let $f\in\Phi$ be given. To see that
\[
FS\big(\{a_n+\sum_{t\in H_n}y_{f(n),t}\}_{n\in\N}\big)\subseteq A,
\]
 let $K\in P_f(\N)$. Let
$K=\{n_1,\cdots,n_m\}$ where $n_1<n_2,\cdots<n_m$. Then $F_{n_1}\subset F_{n_2}\subset\cdots\subset F_{n_m}$
and for each $i\in\{1,\cdots ,m\}$, $\{y_{f(n_i),t}\}_{t\in\N}\in F_{n_i}$ so
\[
\sum_{n\in K}(a_n+\sum_{t\in H_n}y_{f(n),t})=\sum_{i=1}^m(\alpha_\frac{1}{n_i}(F_{n_i})+\sum_{t\in H_\frac{1}{n_i}(F_{n_i})}y_{f(n_i),t})\in A.
\]
The "in particular" is obvious.
\end{proof}
\begin{defn}
Let $S$ be a dense subsemigroup of $\big((0,\infty),+\big)$ and pick $m\in \N$. We define

\begin{multline*}
  \V_m=\{\prod_{i=1}^mH_i\in P_f(\N)^m:\mbox{ if }m>1,\, 1\leq t\leq m-1,\\
  \mbox{ then }maxH_t<minH_{t+1} \}, \\
 \end{multline*}
\[
\J_m=\{\prod_{i=1}^mt(i)\in\N^m:t(1)<\cdots<t(m)\},
\]
and
\[
S^m_\delta=S^m\cap \big(0,\frac{\delta}{m}\big)^m
\]
for $\delta>0$.
\end{defn}
\begin{defn}
 Let $S$ be a dense subsemigroup of $\big((0,\infty),+\big)$. \\
(a) Given $m\in \N$, $\delta>0$, $a\in S_\delta^{m+1}$, $t\in\J_m$, and $f\in\T_0$, define
\[
x(m,a,t,f)=a(m+1)+\sum_{j=1}^ma(j)+f(t(j)).
\]
b) $J_0(S)=\{p\in 0^+: \mbox{ for all } A\in p, A\mbox{ is a }J-\mbox{set near zero}\}.$
\end{defn}
\begin{lem}\label{lemma5.10}
Let $S$ be a dense subsemigroup of $\big((0,\infty),+\big)$ and let $A\subseteq S$. Then $A$ is a $J$-set near zero if and only if
for each $F\in P_f(\T_0)$ and for each $\delta>0$ there exist $m\in \mathbb{N}$, $a\in S_\delta^{m+1}$, and $t\in \J_m$ such that for each $f\in F$, $x(m,a,t,f)\in A$.
\end{lem}
\begin{proof}
 The sufficiency of the statement for $J-$set is trivial. Assume
that $A$ is a $J-$set near zero. Pick $F\in P_f(^\N S)$ and $\delta>0$. Pick $c\in \N$ and for
$f\in F$, define $g_f\in \T_0$ by $g_f(n)=f(n+ c)$. Pick $b\in S$ and
$H\in P_f(\N)$ such that for each $f\in F$, $b+\sum_{t\in H}g_f(t)\in
A$. Let $m=|H|$, and let $t=\big(t(1),\cdots,t(m)\big)$ enumerate $H$ in increasing order, let $a(1)=b$
and for $j\in\{2,\cdots,m+1\}$, let $a(j)=c$. Then the proof is complete.
\end{proof}
\begin{lem}\label{lemma5.11}
Let $S$ be a dense subsemigroup of $\big((0,\infty),+\big)$, $A\subseteq S$ is a $C-$set near zero
if and only if for each $\delta>0$, there exist $m_\delta:P_f(\T_0)\rightarrow \mathbb{N}$,
$\alpha\in\prod_{F\in P_f(\T_0)}S_\delta^{m_\delta(F)+1}$,
and $\tau\in\prod_{F\in P_f(\T_0)}\J_{m_\delta(F)}$ such that \\
(1) if $F,G\in P_f(\T_0)$ and $F\subset G$ then $\tau(F)(m_\delta(F))<\tau(G)(1)$ for each $\delta>0$, and \\
(2) whenever $n\in \mathbb{N}$, $G_1,\cdots, G_n\in P_f(\T_0)$,
$G_1\subset G_2\subset\cdots \subset G_n$, and for each
$i\in\{1,\cdots, n\}$, $f_i\in G_i$, one has
\[
\sum_{i=1}^nx(m_\delta(G_i),\alpha(G_i),\tau(G_i),f_i)\in A.
\]
\end{lem}
\begin{proof}
For the statement about $C-$sets the sufficiency is trivial. For the
necessity, pick $\alpha_\delta:P_f(\T_0)\rightarrow S$ and
$H_\delta:P_F(\T_0)\rightarrow P_f(\N)$ for each $\delta\in (0,1)$
as guaranteed by Definition of
$C-$set near zero. Now pick $c\in \N$ and for $f\in \T_0$ define $g_f\in \T_0$ by
$g_f(s)=f(s+ c)$, for $s\in \N$. For $F\in P_f(\T_0)$ we define
inductively on $|F|$ a set $K(F)\in P_f(\T_0)$ such that\\
(1) $\{g_f:f\in F\}\subseteq K(F)$ and \\
(2) if $\emptyset\neq G\subset F$, then $K(G)\subset K(F)$.

If $F=\{f\}$, let $K(F)=\{g_f\}$. Now let $|F|>1$ and $K(G)$ has been defined for all proper nonempty
subsets of $F$. Pick $h\in \T_0\setminus\bigcup\{K(G):\emptyset\neq G\subset F\}$ and
let $K(F)=\{h\}\cup\{g_f:f\in F\}\cup\bigcup\{K(G):\emptyset\neq G\subset F\}$.

Now for each $\delta\in (0,1)$, we define
  $m_\delta:P_f(\T_0)\rightarrow \mathbb{N}$, $\alpha^\prime\in\prod_{F\in P_f(\T_0)}S^{m_\delta(F)+1}$,
and $\tau\in\prod_{F\in P_f(\T_0)}\J_{m_\delta(F)}$.

Let $F\in P_f(\T_0)$ be given and let $m_\delta(F)=|H_\delta(K(F))|$. Define
$\alpha^\prime (F)\in S^{m_\delta(F)+1}$ by, for $j\in\{1,2,\cdots,
m_\delta(F)+1\}$, $\alpha^\prime (F)(j)=\alpha_\delta(K(F))$ if $j=1$ and
$\alpha^\prime (F)(j)=c$ if $j>1$. Let $\tau(F)=(\tau(F)(1),\cdots,
\tau(F)(m_\delta(F))$ enumerate $H_\delta(K(F))$. We need to show that \\
(1) if $F, G\in P_{f}(\T_0)$ and $F\subset G$, then $\tau(F)(m_\delta(F))<\tau(G)(1)$ for each $\delta\in (0,1)$, and\\
(2) whenever $n\in \mathbb{N}$, $G_1,\cdots, G_n\in P_f(\T_0)$,
$G_1\subset G_2\subset\cdots \subset G_n$, and for each
$i\in\{1,\cdots, n\}$, $f_i\in G_i$, one has
\[
\sum_{i=1}^nx(m_\delta(G_i),\alpha^\prime(G_i),\tau(G_i),f_i)\in A.
\]
 To verify (1), let $F, G\in P_{f}(\T_0)$ with
$F\subset G$, then $K(F)\subset K(F)$, and so $\tau(F)(m_\delta(F))=max H_\delta(K(F))<min H_\delta(K(G))=\tau(G)(1)$.

To verify (2), let $n\in \mathbb{N}$, $G_1,\cdots, G_n\in P_f(\T_0)$,
$G_1\subset G_2\subset\cdots \subset G_n$, and for each
$i\in\{1,\cdots, n\}$, let $f_i\in G_i$. Then $K(G_1)\subset
K(G_2)\subset\cdots \subset K(G_n)$, and for each $f_i\in G_i$,
$g_{f_{i}} \in K(G_{i})$ so
$\sum^{n}_{i=1}(\alpha_\delta(K(G_{i})))+\sum_{t\in H_\delta(K(G_{i}))}g_{f_{i}}(t))\in A$ and\\
\begin{align*}
\sum^{n}_{i=1} & (\alpha_\delta(K(G_{i}))+ \sum_{t \in H_\delta(K(G_{i}))}
g_{f_{i}}(t))
\\
&=\sum^{n}_{i=1}(\alpha_\delta(K(G_{i}))+\sum^{m_\delta(G_{i})}_{j=1}(f_{i}(\tau(G_{i})(j)+c))\\
&=\sum^{n}_{i=1}(\sum^{m_\delta(G_{i})}_{j=1}(\alpha^\prime(G_{i})(j)+(f_{i}(\tau(G_{i})(j)))+\alpha^\prime(G_{i})(m_\delta(G_{i})+1))\\
& =\sum^{n}_{i=1}x(m_\delta(G_i),\alpha^\prime(G_i),\tau(G_i),f_i).
\end{align*}
\end{proof}
\begin{lem}
Let $S$ be a dense subsemigroup of $\big((0,\infty),+\big)$ and let $A\subseteq S$. Let $A$ be a J-set near zero in $S$, then for each
$F\in P_{f}(\T_0)$, each $\delta>0$ and each $n\in\N$, there exist $m\in
\mathbb{N}$, $a\in S_\delta^{m+1}$, and $y\in \J_m $ such that $y(1)>n$ and for each $f\in F$, $x(m,a,y,f)\in A$.
\end{lem}
\begin{proof}
Pick $F\in P_f(\T_0)$, $\delta>0$ and $n\in \N$. For each $f\in F$ define $g_f\in \T_0$ by, for
$u\in \N$, $g_f(u)=f(u+n)$. Pick $m\in \mathbb{N}$, $a\in S_\delta^{m+1}$ and $t\in\J_m$ such that
for each $f\in F$, $x(m,a,t,g_f)\in A$. Define $y\in\J_m$ by $y(i)=n+t(i)$ for $i\in\{1,2,\cdots,m\}$.
Then $y(1)>1$ and for each $f\in F$, $x(m,a,y,f)\in A$.
\end{proof}
\begin{thm}
Let $S$ be a dense subsemigroup of $\big((0,\infty),+\big)$. Then $J_0(S)$ is a compact two sided ideal of
$\beta S$.
\end{thm}
\begin{proof}
Trivially $J_0(S)$ is topologically closed in $\beta S$. Let $p\in
J_0(S)$ and let $q\in \beta S$. We show $q+ p\in J_0(S)$ and $p+
q\in J_0(S)$.\\
To see $q+ p\in J_0(S)$, let $A\in q+p$ and let $F\in P_{f}(\T_0)$.
Then $\{b\in S: -b+A\in p\}\in q$ so pick $b\in S$ such that
$-b+A\in p$. Pick $m\in\mathbb{N}$, $\delta\in (0,1)$, $a\in S_\delta^{m+1}$, and $t\in\J_{m}
$ such that for $f\in F$, $x(m,a,t,f)\in -b + A$. Define $c\in
S_\delta^{m+1}$ by $c(1)=b+a(1)$ and $c(j)=a(j)$ for $j\in
\{2,3,...,m+1\}$. Then for each $f\in F$, $x(m,c,t,f)\in A$.

To see $p+ q\in J_0(S)$, let $A\in p+q$ and let $B=\{x\in S: -x+A\in q \}$.
Then $B\in p$ so for $F\in P_f(\T_0)$ and $\delta\in (0,1)$, pick $m\in\mathbb{N}$, $a\in S_\delta^{m+1}$, and $t\in
\J_{m}$ such that for $f\in F$, $x(m,a,t,f)\in B$. Then
$\bigcap_{f\in F}(-x(m,a,t,f)+A)\in q$ so pick $b\in \bigcap_{f\in
F}-x(m,a,t,f)+A $. Define $c\in S_\delta^{m+1}$ by $c(m+1)= a(m+1)+b$ and
$c(j)=a(j)$ for $j\in \{1,2,...,m\}$. Then for $f\in F $,
$x(m,c,t,f)\in A$.
\end{proof}
\begin{lem}\label{lemma2}
Let $S$ be a dense subsemigroup of $\big((0,\infty),+\big)$. Pick $\delta\in (0,1)$, and let
$m,r\in \mathbb{N}$, let $a\in S_\delta^{m+1}$,
let $t\in \J_{m}$, and for each $y\in \mathbb{N}$, let $c_{y}\in
S_\delta^{r+1}$ and $z_{y}\in \J_{r}$ be a such that for each $y\in
\mathbb{N}$, $z_{y}(r)< z_{y+1}(1)$. Then there exist
$u\in \mathbb{N}$, $d\in S_\delta^{u+1}$, and $q\in \J_{u}$ such that for
each $f\in \T_0$,
\[
(\sum_{j=1}^{m}a(j)+ x(r,c_{t}(j),z_{t}(j),f))+a(m+1)=x(u,d,q,f).
\]
\end{lem}
\begin{proof}
We have
\begin{multline*}
 \sum_{j=1}^{m}(a(j)+x(r,c_{t(j)},z_{t(j)},f))+a(m+1) = \\
 \sum_{j=1}^{m}\big(a(j)+c_{t(j)}(r+1)+\sum_{k=1}^{r}(c_{t(j)}(k)+f(z_{t(j)}(k)))\big)+a(m+1) = \\
 \sum_{j=1}^{m}a(j)+\sum_{j=1}^{m}c_{t(j)}(r+1)+\sum_{j=1}^{m}\sum_{k=1}^{r}(c_{t(j)}(k)+f(z_{t(j)}(k))+a(m+1).\\
\end{multline*}
Now let $u=m\cdot r$. For $j\in \{1,2,...,m\}$ and $p\in
\{1,2,...,r\}$, let $q((j-1)\cdot r+p)=z_{t(j)}(p)$. Let $d(1)=a(1)+
c_{t(1)}(1)$, let $d(u+1)=c_{t(m)}(r+1)+a(m+1)$, for $j\in
\{1,2,...,m-1\}$, let $d(j\cdot
r+1)=c_{t(j)}(r+1)+a(j+1)+c_{t(j+1)}(1)$, and for $j\in
\{1,2,...,m\}$ and $p\in \{2,3,...,r\}$, let $d((j-1)\cdot
r+p)=c_{t(j)}(p)$. So this complete the proof.
\end{proof}
\begin{lem}\label{lemma6}
Let $S$ be a dense subsemigroup of $\big((0,\infty),+\big)$, and let $A_1$ and $A_2$ be subsets of $S$.
 If $A_1 \cup A_2$ is a $J$-set near zero, then either $A_1$ is a J-set near zero or $A_2$ is a J-set neat zero.
\end{lem}
\begin{proof}
Suppose not and pick $F_{1}$ and $F_{2}$ in $P_{f}(^\N S)$ and $\delta>0$ such that
for each $i\in \{1,2\}$, each $u \in \mathbb{N}$, each $d\in S_\delta^
{u+1}$, and each $q\in \J_{u}$, there is some $f\in F_{i}$
such that $x(u,d,q,f) \notin A_{i}$.

Let $F=F_{1}\cup F_{2}$, $k=\mid F\mid$, and write
$F=\{f_{1},f_{2},...,f_{k}\}$. Pick by Lemma, some $n\in \mathbb{N}$
such that whenever length $n$ words over the alphabet
$\{1,2,...,k\}$ are 2-colored, there is a variable word $w(v)$
beginning and ending with a constant and without successive
occurrences of $v$ such that $\{w(l): l\in \{1,2,...,k\}\}$ is
monochromatic.

Let $W$ be the set of length $n$ words over $\{1,2,...,k\}$. For
$w=b_{1},b_{2},...,b_{n}\\ \in W$(where each $b_{i}\in
\{1,2,...,k\}$), define $g_{w}: \N \rightarrow S$ by, $y\in S$,
$g_{w}(y)=\sum^{n}_{i=1} f_{b_{i}}(ny+ix)$ where $x\in S$. Since $A$
is a $J_S$-set near zero, pick $m\in \mathbb{N}$, $a\in S_\delta^{m+1}$,
 $t \in \J_m$
 such that for all $w\in W$, $x(m,a,t,g_{w})\in A$. Define
$\varphi: W \rightarrow \{1,2\}$ by $\varphi (w)=1$ if
$x(m,a,t,g_{w})\in A_{1}$ and $\varphi (w)=2$ otherwise. Pick a
variable word $w(v)$, beginning and ending with a constant and
without successive occurrences of $v$ such that $\varphi$ is
constant on $\{w(l):\l\in \{1,2,...,k\}\}$. Assume without loss of
generality that $\varphi (w(l))=1$ for all $j\in \{1,2,...,k\}$.
That is, for all $l\in \{1,2,...,k\}$,
\[
(\sum ^{m}_{j=1}a(j) +
g_{w(l)}(t(j))) + a(m+1)=x(m,a,t,g_{w(l)})\in A_{1}.
\]

Let $w(v)=b_{1}b_{2}\cdots b_{n}$ where each$b_{i}\in
\{1,2,,...,k\}\cup \{v\}$, some $b_{i}=v$, $b_{1}\neq v$, $b_{n}\neq
v$, and if $b_{i}=v$, then $b_{i+1}\neq v$. Let $r$ be the number of
occurrences of $v$ in $w(v)$ and pick $L\in \V_{r+1}$ and $s\in\J_r$
such that for each $p\in\{1,\cdots,r\}$, $max L_{p}< s(p)<min L_{p+1}$,
\[
\bigcup_{p=1}^{r+1}L_p=\{i\in\{1,\cdots,n\}:b_i\in\{1,\cdots,k\}\}
\]
and $\{s(1),\cdots,s(r)\}=\{i\in\{1,\cdots,n\}:b_i=v\}$. (For example, if $w(v)=12v131v2v1121v32$,
then $r=4$, $L=(\{1,2\},\{4,5,6\},\{8\},\{10,11,12\},\{14,15\})$,
and
$s=(3,7,9,14)$.)

We shall show now that, given $y\in \N$, there exist $c_{y}\in
S_\delta^{r+1}$ and $z_{y}\in \J_{r}$ such that for all $l\in
\{1,2,...,k\}$, $g_{w(l)}(y)=x(r,c_{y},z_{y},f_{l})$ and further,
for each $y$, $z_{y}(r)<z_{y+1}(1)$. So let $y\in \N$ be given. For $p\in
\{1,2,\cdots,r+1\}$, let $c_{y}(p)=\sum _{i\in L(p)}f_{b_{i}}(ny+i)$
and for $p\in \{1,2,\cdots,r\}$, let $z_{y}(p)=ny+s(p)$. To see that
these are as required, first note that  $z_{y}(r)\leq ny+n<z_{y+1}(1)$. Now let $l\in\{1,\cdots,k\}$ be given.
Then $w(l)=d_{1}d_{2}\cdots d_{n}$ where for $i\in \{1,2,\cdots,n\}$, $d_{i}=b_{i}$ if $i\in \bigcup^{r+1}_{p=1}L_p$ and $d_{i}=l$
if $i\in \{s(1),s(2),\cdots,s(r)\}$.

Therefore
\begin{align*}
g_{w(l)}(y) &=\sum^{n}_{i=1}f_{d_{i}}(ny+i) \\
 &=(\sum^{r}_{p=1}(\sum _{i\in L(p)}f_{b_{i}}(ny+i))+ f_{l}(ny+s(p)))+
\sum_{i\in L(r+1)}f_{b_{i}}(ny+i)\\
 &=(\sum^{r}_{p=1}c_{y}(p)+ f_{l}(z_{y}(p)))+ c_{y}(r+1)\\
&=x(r,c_{y},z_{y},f_{l})
\end{align*}
as required.

Now pick Lemma \ref{lemma2}, $u\in\N$, $d\in S_\delta^{u+1}$, $q\in
\mathcal{J}_{u}$ such that for each $f\in ^\N S$,
\[
(\sum^{m}_{j=1}a(j)+ x(r,c_{t(j)},z_{t(j)},f)) +
a(m+1)=x(u,d,q,f).
\]
Pick $l\in \{1,2,\cdots,k\}$ such that $f_{l}\in
F_{1}$ and $x(u,d,q,f_{l})\notin
A_{1}$. But
\begin{center}
\begin{align*}
x(u,d,q,f_{l}) &=(\sum ^{m}_{j=1}a(j)+
x(r,c_{t(j)},z_{t(j)},f_{l}))+ a(m+1) \\
&=(\sum^{m}_{j=1}a(j)+ g_{w(l)}(t(j)))+ a(m+1)\\
&=x(m,a,t,g_{w(l)})\in A_{1},
\end{align*}
\end{center}
a contradiction.
\end{proof}
\begin{thm}\label{theorem6}
Let $S$ be a dense subsemigroup of $\big((0,\infty),+\big)$, let $A\subseteq S$. Then $\overline{A}\cap
J_0(S)\neq \emptyset$ if and only if $A$ is a $J_S$-set near zero.
\end{thm}
\begin{proof}
The necessity is trivial. By Lemma \ref{lemma6}, $J_S$-sets are partition regular.
So, if $A$ is a $J_S$-set near zero, by Theorem \ref{theorem0}, there is some $p\in \beta S$
such that $A\in p$ and for every $B\in p$, $B$ is a $J_S$-set near zero.
\end{proof}
\begin{cor}
Let $S$ be a dense subsemigroup of $\big((0,\infty),+\big)$, and let $A$ be a piecewise syndetic near zero
subset of $S$. Then $A$ is a $J_S$-set near zero.
\end{cor}
\begin{proof}
By Theorem \ref{theorem01}, $\overline{A}\cap K(0^+(S))\neq \emptyset$. Since $K(0^+(S))\subseteq J_0(S)$,
so $\overline{A} \cap J_0(S) \neq
\emptyset$ so by Theorem \ref{theorem6}, $A$ is a $J_S$-set near zero.
\end{proof}

\begin{thm}
Let $S$ be a semigroup and let $A\subseteq S$. If there is an idempotent in $\overline{A}\cap
J_0(S)$, then $A$ is a C-set near zero.
\end{thm}
\begin{proof}
Pick $p=p+p\in \overline{A}\cap J_0(S)$. Recall that $A^*=\{x\in A:-x+A\in p\}$ and, by lemma \ref{lemma11}, if $x\in A^*$,
then $-x+A^*\in p$. For  every $\delta \in (0,1)$ we define $m_\delta(F)$ and $\alpha(F)$ and $\tau(F)$
for $F\in P_f(\T_0)$ by induction on $|F|$  so that

(1) if $F,G\in P_f(\T_0)$ and $F\subset G$ then $\tau(F)(m_\delta(F))<\tau(G)(1)$ for each $\delta>0$, and

(2) whenever $n\in \mathbb{N}$, $G_1,\cdots, G_n\in P_f(\T_0)$,
$G_1\subset G_2\subset\cdots \subset G_n$, and for each
$i\in\{1,\cdots, n\}$, $f_i\in G_i$, one has
\[
\sum_{i=1}^nx(m_\delta(G_i),\alpha(G_i),\tau(G_i),f_i)\in A.
\]
Assume  first that $F=\{f\}$. Then $A^*$ is a $J_S$-set near zero so pick $m_\delta(F)\in  \mathbb{N}$, $\alpha(F)\in S^{m_\delta(F)+1}$,
and $\tau(F)\in \mathcal{J}_{m_\delta(F)}$ such that
\[
x(m_\delta(F),\alpha(F),\tau(F),f)\in A^*.
\]
 Now assume that $|F|>1$
and that $m_\delta(G)$, $\alpha(G)$, and $\tau(G)$ have been defined for all non-empty proper subsets $G$ of $F$ and for each $\delta>0$.
Pick $\delta>0$, and let $k=max\{\tau(G)(m_\delta(G)): \emptyset\neq G\subsetneq F\}$.
Let
\begin{multline*}
 M_\delta =\{\sum_{i=1}^m x(m_\delta(G_i),\alpha(G_i),\tau(G_i),f_i):n\in \mathbb{N},\\
\emptyset \neq G_1\subset G_2\subset\cdots\subset G_n=F,
and\,\, \{f_i\}_{i=1}^n\in\Pi_{i=1}^n G_i\}.\\
\end{multline*}
Let $B=A^*\cap\bigcap_{b\in M_\delta}(-b+A^*) $. Since $M_\delta$ is a finite subset of $A^*$, $B\in p$
 and therefor $B$ is a J-set near zero. Pick by Lemma \ref{lemma5.10}, $m_\delta(F)\in \mathbb{N}$, $\alpha(F)\in S^{m_\delta(F)+1}$ ,
 and $\tau(F)\in \mathcal{J}_{m_\delta(F)}$ such that $\tau(F)(1)>k$ and for each $f\in F$, $x(m_\delta(F),\alpha(F),\tau(F),f)\in B$.

 Hypothesis $(1)$ is satisfied directly. To verify hypothesis $(2)$, let $n\in \mathbb{N}$, let $\emptyset \neq G_1\subset G_2\subset\cdots\subset G_n=F$,
 and for each $i\in \{1,2,\cdots,n\}$, let $f_i\in G_i$. If $n=1$, then $x(m_\delta(G_1),\alpha(G_1),\tau(G_1),f_i)\in B\subseteq A^*$,
 so assume that $n>1$. Let $b=\sum_{i=1}^{n-1} x(m_\delta(G_i),\alpha(G_i),\tau(G_i),f_i)$ then $b\in M_\delta$ so
 $x(m_\delta(G_n),\alpha(G_n),\tau(G_n),f_i)\in B\subseteq -b+A^*$ so
 \[
 \sum_{i=1}^{n-1} x(m_\delta(G_i),\alpha(G_i),\tau(G_i),f_i)\in A^*
 \]
  as required.
\end{proof}
\begin{cor}
Let $S$ be a dense subsemigroup of $((0,\infty),+)$ and let $A$ be a central set near zero in $S$. Then $A$ is
a $C$-set near zero.
\end{cor}
\begin{proof}
It is obvious.
\end{proof}
\begin{thm}
Let $S$ be a dense subsemigroup of $((0,\infty),+)$, let $A$ be a central subset near zero of $S$, and for each
$l\in\N$, let $\{y_{l,n}\}_{n\in\N}$ be a sequences in $S$
such that $lim_{n\rightarrow \infty}y_{l,n}=0$ for each $l\in\N$. Given $l,m\in\N$,
$a\in S^{m+1}_{\frac{1}{m}}$, and $H\in\V_m$, let
\[
w(a,H,l)=\big(\sum_{i=1}^m\big(a(i)+\sum_{t\in H(i)}y_{l,t}\big)\big)+a(m+1).
\]
There exist sequences $\{m(n)\}_{n\in\N}$, $\{a_n\}_{n\in\N}$, and $\{H_n\}_{n\in\N}$ such that

(1) for each $n\in\N$, $m(n)\in\N$, $a_n\in S_{\frac{1}{m(n)}}^{m(n)+1}$, $H_n\in\V_{m(n)}$, and
$maxH_{n,m(n)}<minH_{n+1,1}$, and

(2) for each $f\in\Phi$, $FS\big(\{w(a_n,H_n,f(n))\}_{n\in\N}\big)\subseteq A$.
\end{thm}
\begin{proof}
As in the proof of Theorem \ref{theorem5.7}, we may assume that the sequences $\{y_{l,n}\}_{n\in\N}$ are
all distinct. $A$ is central near zero, so $A$ is a $C$-set. For each $k\in\N$, pick $m_\frac{1}{k}^\prime:P_f(\T_0)\rightarrow \mathbb{N}$,
$\alpha\in\prod_{F\in P_f(\T_0)}S_\frac{1}{k}^{m_\frac{1}{k}^\prime(F)+1}$,
and $\tau\in\prod_{F\in P_f(\T_0)}\J_{m_k^\prime(F)}$ as guaranteed by the fact that $A$ is a $C$-set. For each $n\in\N$, let
$F_n=\{\{y_{1,t}\}_{t\in\N},\cdots,\{y_{n,t}\}_{t\in\N}\}$, $m(n)=m^\prime_\frac{1}{n}(F_n)$, $a_n=\alpha(F_N)$ and
\[
H_n=\big(\{\tau(F_n)(1)\},\cdots,\{\tau(F_n)(m(n))\}\big).
\]
To see that $m(n)$, $a_n$ and $H_n$ are as required, let $f\in\Phi$ and let $K\in P_f(\N)$.
Enumerate $K$ in order as $n(1),n(2),\cdots,n(l)$. For $i\in\{1,\cdots,l\}$, let $g_i=\{y_{f(n(i)),t}\}_{t\in\N}$. Then
\[
\sum_{n\in K}w(a_n,H_n,f(n))=\sum_{i=1}^lx\big(m^\prime_\frac{1}{n(i)}(F_{n(i)}),\alpha(F_{n(i)}),\tau(F_{n(i)}),g_i\big)\in A.
\]
\end{proof}
\begin{lem}\label{lemma3.21}
Let $J$ be a set, let $(D,\leq)$ be a directed set, and let $S$ be a dense subsemigroup of $((0,\infty),+)$.
Let $\{T_i\}_{i\in D}$ be a decreasing family of nonempty subsets of $S$ such that

1) $0\in cl_{\mathbb{R}}T_i$,

2) $\bigcap_{i\in D}T_i=\emptyset$, and

3) for each $i\in D$ and each $x\in T_i$ there is some $j\in D$ such that $x+T_j\subseteq T_i$.

 Let
$Q=\bigcap_{i\in D}cl_{\beta S_d}T_i$. Then $Q$ is a compact subsemigroup of $0^+(S)$. Let $\{E_i\}_{i\in D}$ and
$\{I_i\}_{i\in D}$ be decreasing families of nonempty subsets of $\Pi_{t\in J}S$ with the following
properties:

$(a)$ for each $i\in D$, $I_i\subseteq E_i\subseteq \Pi_{t\in J} T_i$,

$(b)$ for each $i\in D$ and each $\overrightarrow{x} \in I_i$ there exists $j\in D$ such that $\overrightarrow{x}+E_j\subseteq I_i$,
and

$(c)$ for each $i \in D$ and each $\overrightarrow{x}\in E_i\setminus I_i$ there exists $j\in D$ such that
$\overrightarrow{x}+ E_j\subseteq E_i$ and $\overrightarrow{x} + I_j\subseteq I_i$.

Let $Y=\Pi_{t\in J}0^+(S)$, let $E=\bigcap_{i\in D}cl_YE_i$, and let $I=\bigcap_{i\in D}cl_Y I_i$. Then $E$ is
a subsemigroup of $\Pi_{t\in J}Q$ and $I$ is an ideal of $E$. If, in addition, either

$(d)$ for each $i \in D$, $T_i=S$ and $\{a\in S:\overline{a}\notin E_i\}$ is not piecewise syndetic near zero, or

$(e)$ for each $i\in D$ and each $a\in T_i$ , $\overline{a}\in E_i$,

then given
any $p\in K(Q)$, one has $\overline{p}\in E\cap K(\Pi_{t\in J}Q)=K(E)\subseteq I$.
\end{lem}
\begin{proof}
By Theorem 4.20 in \cite{hinbook}, $Q$ is a subsemigroup of $0^+(S)$. For the proof that $E$ is
a subsemigroup of $\Pi_{t\in J}Q$ and $I$ is an ideal of $E$, see the proof of Lemma 14.9 in \cite{hinbook}.

To complete the proof, assume that $(d)$ or $(e)$ holds. It suffices to establish
\[
\mbox{if }p\in K(Q),\mbox{ then }\overline{p}\in E.\,\,\,\,\,(*)
\]
Indeed, assume we have established $(*)$. Then $\overline{p}\in E\cap\Pi_{t\in J}K(Q)$ and
$\Pi_{t\in J}K(Q)=K(\Pi_{t\in J}Q)$ by Theorem 2.23 in \cite{hinbook}. Then by Theorem 1.65 in \cite{hinbook},
$K(E)=E\cap K(\Pi_{t\in J}Q)$ and, since $I$ is an ideal of $E$, $K(E)\subseteq I$.

To establish $(*)$, let $p\in K(Q)$ be given. To see that $\overline{p}\in E$, let $i\in D$ be given and
let $U$ be a neighborhood of $\overline{p}$. Pick $F\in P_f(J)$ and for each $t\in F$ pick
some $A_t\in p$ such that $\bigcap_{t\in F}\pi_t^{-1}[cl_{\beta S_d}A_t]\subseteq U$, where $\pi_t$ is projection for $t\in J$.

Assume now that $(d)$ holds. Since $p\in K(0^+(S)$ and $\{a\in S:\overline{a}\notin E_i\}$ is not piecewise syndetic near zero,
so by Theorem 2.4, $\{a\in S:\overline{a}\notin E_i\}\notin p$ and hence $\{a\in S:\overline{a}\in E_i\}\in p$. Now
pick $a\in\big(\bigcap _{t\in F}A_t\big)\cap\{a\in S:\overline{a}\in E_i\}$. Then $\overline{a}\in U\cap E_i$.

If $(e)$ holds, see the proof of Lemma 14.9 in \cite{hinbook}.
\end{proof}
\begin{thm}
Let $S$ be a dense subsemigroup of $((0,\infty),+)$ and let $A\subseteq S$. Then $A$ is a C-set near zero
if and only if there is an idempotent in $\overline{A}\cap J_0(S)$.
\end{thm}
\begin{proof}
The sufficiency is obvious.

Pick for every $k\in\N$,
there exist $m_\frac{1}{k}:P_f(\T_0) \rightarrow \mathbb{N}$, $\alpha\in \Pi _{F\in P_f(\T_0)}S^{m_\frac{1}{k}(F)+1}$ ,
and $\tau\in \Pi_{F\in P_f(\T_0)}\mathcal{J}_{m_\frac{1}{k}(F)}$ as for each $\frac{1}{k}\in (0,1]$ guaranteed by the
fact that $A$ is a C-set near zero.
For $F\in P_f(\T_0)$ and $k\in\N$ define
\begin{multline*}
T_{F,k}=\{\sum_{i=1}^n x(m_\frac{1}{k}(F_i),\alpha(F_i),\tau(F_i),f_i):n\in \mathbb{N},
\forall F_i\in P_f(\T_0),\\
  F_1\subset F_2\subset \cdots\subset F_n,\mbox{ and  for  each }\\
  i\in \{1,2,...,\},  f_i\in F_i\}.\\
   \end{multline*}
It is obvious that if $F,G\in P_f(\T_0)$, then $T_{F\cup G,k}\subseteq T_{F,k} \cap T_{G,k}$ for each $k\in\N$,
so $Q_k=\bigcap_{F\in P_f(\T_0)}\overline{T_{F,k}}$ is a non-empty set for each $k\in\N$.
Pick $k\in\N$, we show $Q_k$ is a subsemigroup of $0^+(S)$.

For this it suffices by Theorem 4.20 in \cite{hinbook}, we show that for all $F\in P_f(\T_0)$ and
all $u\in T_{F,k}$, there is some $G\in P_f(\T_0)$ such
that $u+T_{G,k}\subseteq T_{F,k}$. So let $F\in P_f(\T_0)$ and $u\in T_F$ be
given. Pick $n\in \mathbb{N}$, strictly increasing $\{F_i\}_{i=1}^{n}$ in $P_f(\T_0)$ such that
$F\subset F_1$, and $f\in \Pi_ {i=1}^{n}F_i$ such that
\[
u=\sum_{i=1}^n x(m_\delta(F_i),\alpha(F_i),H_\delta(F_i),f_i).
\]
Then $u+T_{F_n,k}\subseteq T_{F,k}$. Therefore for each $k\in\N$, $Q_k$ is a compact subsemigroup of $0^+(S)$.

Now for some $k\in \N$, we show that $K(Q_k)\subseteq \overline{A}\cap J_0(S)$ so that any idempotent in $K(Q)$ establishes
the theorem. We have that each $T_{F,k}\subseteq \overline{A}$ so $Q_k\subseteq \overline{A}$. Let $p\in K(Q_k)$.
 We need to show that $p\in J_{0}(S)$, so let $B\in p$. We shall show that $B$ is a $J$-set near zero.
 So let $F\in P_f(\T_0)$. We shall produce $v\in \mathbb{N}$
$c\in S_\frac{1}{k}^{v+1}$, and $t\in \mathcal{J}_{v}$ such that for each $f\in F$,  $x(v,c,t,f)\in B$.

We apply Lemma \ref{lemma3.21} with $J=F$ and $D=\{G\in P_f(\T_0) : F\subseteq G\} $. Pick $k\in \N$, and note that
$Q_k=\bigcap_{G\in D}\overline{T_{G,k}}$ as in Lemma \ref{lemma3.21}.  For $G\in D$ we shall define a subset
$I_G$ of $\Pi_{f\in F}T_{G,k}$ as follows. Let $w\in \Pi_{f\in F}T_{G,k}$ then $w\in I_G$ if and only if there is
some $n\in \mathbb{N}-\{1\}$ such that there exist $C_1$, $C_2$, $\{G_i\}_{i=1}^{n}$ and $\eta$ such that

$(1)$ $C_1$ and $C_2$ are disjoint nonempty sets and $\{1,2,...,n\}=C_1\cup C_2$,

$(2)$ $\{G_i\}_{i=1}^{n}$ is strictly increasing in $P_f(\T_0)$ with $G\subset G_1$, and

$(3)$ $\eta\in \Pi_{i\in C_1}G_i$ and for each $f\in F$, if $\gamma_f\in \Pi_{i=1}^{n}G_i$ is defined
by
\begin{center}
\begin{equation*} \gamma_f(i) = \left\{
\begin{array}{ll}
\eta_i & \text{if } i\in C_1\\
f & \text{if } i\in C_2
\end{array} \right.
\end{equation*}
\end{center}
then $w(f)=\sum_{i=1}^{n}x(m_\frac{1}{k}(G_i),\alpha(G_i),\tau(G_i),\gamma_f(i)$.

For $G\in D$, note that $I_G\neq \emptyset$ and let $E_G=I_G\bigcup \{\overline{b}: b\in T_G\}$.

We claim that $\{E_G\}_{G\in D}$ and $\{I_G\}_{G\in D}$ satisfy
statements $(a),(b),(c)$ and $(e)$ of Lemma \ref{lemma3.21}. Statements $(a)$ and $(e)$ hold trivially.

To verify $(b)$, let $G\in D$ and $w\in I_G$. Pick $n,C_1,C_2, \{G_i\}_{i=1}^{n}$ and $\eta$ as
guaranteed by the fact that $w\in I_G$. We claim that $w+E_{G_n}\subseteq I_G$. So let $z\in E_{G_n}$.

Assume first that $z=\overline{b}$ for some $b\in T_{G_n}$. Pick $n^\prime \in \mathbb{N}$,
strictly increasing $\{F_i\}_{i=1}^{n^\prime}$ in $P_f(\T_0)$ with $G_n\subset F_1$,
and $\eta^\prime \in \Pi_{i=1}^{n^\prime}F_i$ such that
\[
b=\sum_{i=1}^n x(m_\delta(F_i),\alpha(F_i),\tau_(F_i),\tau(F_i),\eta^\prime(i)).
\]
Let $C_1^{\prime\prime}=C_1\cup \{n+1,n+2,...,n+n^\prime\}$ and for $i\in \{1,2,...,n+n^\prime\}$ and
\begin{center}
\begin{equation*} L_{i} = \left\{
\begin{array}{ll}
G_{i} & \text{if } i\leq n\\
F_{i-n} & \text{if } i>n.
\end{array} \right.
\end{equation*}
\end{center}
Define $\eta^{\prime\prime} \in \Pi_{i\in G_1^{\prime\prime}}L_i $ by, for $i\in C_1^{\prime\prime}$,
\begin{center}
\begin{equation*}\eta^{\prime\prime}(i) = \left\{
\begin{array}{ll}
\eta(i) & \text{if } i\leq n\\
\eta^{\prime}(i-n) & \text{if } i>n.
\end{array} \right.
\end{equation*}
\end{center}
Then $n+n^\prime$ , $C_1^{\prime\prime},C_2^{\prime\prime},\{G_i\}_{i=1}^{n+n^{\prime}} $, and $\eta^{\prime\prime}$ establish that $w+z\in I_G$.

Now assume that $z\in I_{G_n}$. Pick $n^\prime, C_1^\prime, C_2^\prime,\{F_i\}_{i=1}^{n^\prime}$ and $\eta^\prime$ as guaranteed
by the fact that $z\in I_{G_n}$. Let $C_1^{\prime\prime}=C_1\cup\{n+i:i\in C_1^\prime\}$, let
$C_2^{\prime\prime}=C_2\cup\{n+i:i\in C_2^\prime\}$, and for $i\in\{1,2,\cdots,n+n^\prime\}$ let
\begin{center}
\begin{equation*} L_{i} = \left\{
\begin{array}{ll}
G_{i} & \text{if } i\leq n\\
F_{i-n} & \text{if } i>n.
\end{array} \right.
\end{equation*}
\end{center}
Define $\eta^{\prime\prime}\in\Pi_{i\in C_1^{\prime\prime}}L_i$ by, for $i\in C_1^{\prime\prime}$,
\begin{center}
\begin{equation*} \eta^{\prime\prime}(i) = \left\{
\begin{array}{ll}
\eta(i) & \text{if } i\leq n\\
\eta^\prime(i-n) & \text{if } i>n.
\end{array} \right.
\end{equation*}
\end{center}
Then $n+n^\prime, C_1^{\prime\prime}, C_2^{\prime\prime}, \{L_i\}_{i=1}^{n+n^\prime\prime}$,
and $\eta^{\prime\prime}$ establish that $w+z\in I_G$.

To verify $(c)$ let $G\in D$ and let $w\in E_G \setminus I_G$. pick $b\in T_G$ such that $w=\bar{b}$. Pick $n\in \mathbb{N}$,
strictly increasing $\{G_i\}_{i=1}^{n}$ in $P_f(\T_0)$
with $G\subset G_1$, and $\eta \in \Pi_{i=1}^{n}G_i$ such that $b=\sum_{i=1}^{n}x(m_\frac{1}{k}(G_i),\alpha(G_i),\tau(G_i),\eta(i))$.
Then as above one has that $w+E_{G_n}\subseteq E_G$ and $w+I_{G_n}\subseteq I_G$.

We then have by Lemma \ref{lemma3.21} that $\overline{p}\in \bigcap_{G\in D}\overline{I_G}$.
Now $\Pi_{f\in F}\overline{B}$ is a neighborhood of $\overline{p}$
so pick $w\in I_F \cap \Pi_{f\in F}\overline{B} $. Pick $n,C_1,C_2,\{G_i\}_{i=1}^{n}$
and $\eta\in \Pi_{i\in C_1}G_i$ as guaranteed by the fact that $w\in I_F$. Let $r=|C_2|$ and let $h_1,h_2,\cdots,h_r$
be the elements of $C_2$ listed in increasing order. Let $v=\sum_{i=1}^{r}m_\frac{1}{k}(G_{h_i}) $. If $h_1=1$, let $c(1)=\alpha(G_1)(1)$.
If $h_1>1$, let
\[
c(1)=\sum_{i=1}^{h_{1}-1}x(m_\frac{1}{k}(G_i),\alpha(G_i),\tau(G_i),\eta(i)+\alpha(G_{h_1})(1).
\]
For $1<j\leq m_\delta(G_{h_1})$ let $c(j)=\alpha(G_{h_1})(j)$ and for $1\leq j\leq m_\frac{1}{k}(G_{h_1})$
let $\mu (j)=\tau(G{h_1})(j)$.

Now let $s\in \{1,2,...,r-1\}$ and let $u=\sum _{i=1}^{s}m_\frac{1}{k}(G_{h_i})$ if $h_{s+1}=h_s+1$
let $c(u+1)=\alpha(G_{h_s})(m_\frac{1}{k}(G_{h_s})+\alpha(G_{h_{s+1}})(1)$. If $h_{s+1}>h_s+1$,
let
\begin{multline*}
c(u+1)=\alpha(G_{h_s})(m_\frac{1}{k}(G_{h_s}+1)\\
+(\sum_{i=h_s+1}^{h_{s+1}-1})x(m_\frac{1}{k}(G_i),\alpha(G_i),\tau(G_i),\eta(i))+\alpha (G_{h_{s+1}}(1)).\\
\end{multline*}
And for $u<j\leq \sum _{i=1}^{s+1}m_\frac{1}{k}(G_{h_i})$, let $\mu(j)=\tau(G_{h_{s+1}})(j-u)$.

If $h_r=n$, let $c(v+1)=\alpha(G_n)(m_\frac{1}{k}({G_n+1}))$ if $h_r<n,$ let
 \[
c(v+1)=\alpha(G_n)(m_\frac{1}{k}({G_n+1}))+\sum_{i=h_r}^{n}(x(m_\frac{1}{k}(G_i),\alpha(G_i),\tau(G_i),\eta(i))).
\]
Then $c\in S_\frac{1}{k}^{v+1}$, $M\in \J_v$ and for each $f\in F$, $x(v,c,\mu,f)\in B$ as required.
\end{proof}


\providecommand{\bysame}{\leavevmode\hbox
to3em{\hrulefill}\thinspace}


\end{document}